\documentclass[british,aip,reprint]{ieeeaccess}
\usepackage[T1]{fontenc}
\usepackage[latin9]{inputenc}
\usepackage{array}
\usepackage{multirow}
\usepackage{amsmath}
\usepackage{amsthm}
\usepackage{graphicx}

\makeatletter

\providecommand{\tabularnewline}{\\}

\usepackage{widetext}
\usepackage{colortbl}

\makeatother

\usepackage{babel}
\begin{document}
\emph{\small{}\history{Date of publication xxxx 00, 0000, date of current version xxxx 00, 0000.} \doi{10.1109/ACCESS.2017.DOI}
\title{Dynamics of targeted ransomware negotiation} 
\author{\uppercase{Pierce Ryan}\authorrefmark{1,2}, \IEEEmembership{Graduate Student Member, IEEE}, \uppercase{\uppercase{John Fokker}\authorrefmark{3}, \uppercase{Sorcha Healy}\authorrefmark{4}, and \uppercase{Andreas Amann}\authorrefmark{1}}} 
\address[1]{School of Mathematical Sciences, University College Cork, Cork T12 XF62, Ireland} 
\address[2]{Artificial Intelligence Research, McAfee, T12 RRC9 Cork, Ireland.} 
\address[3]{Trellix Labs, Trellix, 1119PE, Schiphol-Rijk, Netherlands.} 
\address[4]{Global Data Science \& Analytics, Microsoft Ireland, Dublin D18 P521, Ireland.} 
\tfootnote{This work was supported by the Irish Research Council and McAfee LLC through the Irish Research Council Employment-Based Postgraduate Programme under grant EBPPG/2018/269.}
\markboth {Author \headeretal: Preparation of Papers for IEEE TRANSACTIONS and JOURNALS} {Author \headeretal: Preparation of Papers for IEEE TRANSACTIONS and JOURNALS}
\corresp{Corresponding author: Pierce Ryan (e-mail: pierceryan@umail.ucc.ie).}
\begin{abstract} In this paper, we consider how the development of targeted ransomware has affected the dynamics of ransomware negotiations to better understand how to respond to ransomware attacks. We construct a model of ransomware negotiations as an asymmetric non-cooperative two-player game. In particular, our model considers the investments that a malicious actor must make in order to conduct a successful targeted ransomware attack. We demonstrate how imperfect information is a crucial feature for replicating observed real-world behaviour. Furthermore, we present optimal strategies for both the malicious actor and the target, and demonstrate how imperfect information results in a non-trivial optimal strategy for the malicious actor. \end{abstract}
\titlepgskip=-15pt
\maketitle}{\small\par}

\section{Introduction}

Computer security is a rapidly developing field, with new threats
emerging and evolving constantly. As computer security providers develop
their methods for detecting malware (malicious software), the malicious
actors behind the various strains of malware are forced to refine
their techniques for avoiding detection, prompting further development
from the computer security industry. As a result of the interaction
between these competing agendas, problems in computer security can
give rise to rich dynamical behaviour. By analysing these dynamics,
we can provide insights that assist in understanding phenomena observed
in computer security. The current paper seeks to provide insight on
recent developments in ransomware by using game theory to explore
the dynamics they have introduced.

Ransomware is a type of malware designed to extort a ransom from the
victim \cite{kalaimannan2017influences,maigida2019systematic}, usually
by denying the victim access to their computer or data until the ransom
has been paid. In the past, ransomware relied on extorting a small
amount of money from a large number of victims. The ransom itself
would be fixed at a price low enough that nearly anyone could pay,
but was typically non-negotiable, as negotiating a small ransom with
many victims would not be worth the effort. In recent years, a new
phenomenon known as \emph{targeted ransomware} has emerged \cite{beek2016targeted,bajpai2020dissecting}.
Malicious actors operating targeted ransomware specifically target
large organisations, which can be extorted for significantly higher
ransoms \cite{coveware2020ransomware,zimba2019economic}. However,
they are also likely to have a higher level of computer security,
and so the malicious actors are forced to make significant investment
into breaching their security. Given the effort involved in breaching
security, the malicious actors will invest further in calculating
the highest ransom they believe their target will pay. As a result
of the larger sums of money at stake, malicious actors are willing
to negotiate the ransom demand to facilitate payment. We consider
these negotiations to be a crucial feature in targeted ransomware,
as their outcome has an enormous impact on both the malicious actors
and the targeted organisation. In this paper, we develop a model of
targeted ransomware negotiations based on game theory. Game theory
is a branch of mathematics which studies strategic interactions between
rational decision-makers \cite{myerson2013game}. A \emph{game} is
a mathematical model with a clearly defined set of rules where two
or more \emph{players} make strategic decisions to influence the outcome
of the game to their own personal benefit. By analysing the decisions
available to players, optimal decision-making strategies can be determined
which offer the best outcome for each player. Game theory can be applied
to many real-life scenarios that involve competing interests by formulating
a game as a mathematical abstraction of the given scenario. By analysing
the decisions taken in the game logically, optimal strategies can
be determined, granting insight into real-world behaviour. As such,
game theory is highly applicable to fields such as ecology \cite{brown1999ecology,mcgill2007evolutionary,smith1973logic},
economics \cite{friedman1998economic,selten1990bounded,lukas2012earnouts,cerdeiro2017contagion}
and politics \cite{kydd1997game,ward1993game}. Very recently, ideas
from game theory have been found to be useful in the study of ransomware.
While all ransomware follows the same fundamental principles, there
is sufficient variety in observed phenomena to merit a variety of
models, such as defence and deterrence \cite{caulfield2015optimizing,lindsay2015tipping,laszka2017economics,cartwright2019pay,hu2020optimal},
iterative negotiations \cite{caporusso2018game}, price discrimination
\cite{hernandez2020economic}, incentive to return encrypted data
data \cite{cartwright2019ransomware}, and sale of stolen data \cite{li2021game}.
In this paper, we examine how game theory can be applied to the growing
threat of targeted ransomware. We construct a model of ransom negotiations
as a game played between a malicious actor and their target. The game
focuses on the strategic behaviour and investments that the malicious
actor must commit to in order to implement a successful targeted ransomware
attack. By analysing our model, we demonstrate how imperfect information
is crucial for replicating observed real-world behaviour, and provide
new insights into the real-world behaviour and strategy of malicious
actors operating strains of targeted ransomware. 

\section{Background}

In order to construct a model of targeted ransomware negotiations,
we first consider the key features of targeted ransomware. On a technical
level, one of the main differences between untargeted and targeted
ransomware is how it spreads \cite{panda2019targeted}. Untargeted
ransomware is distributed indiscriminately, relying on victims with
weak security to make a mistake that allows the ransomware to infect
their computer. Such a strategy may be appealing to a malicious actor,
as it requires a low investment of effort to infect victims. However,
this indiscriminate strategy is relatively easy to defend against.
A potential victim can reduce their chances of being infected through
practices such as maintaining good security, careful internet usage,
and maintaining up-to-date data backups. Large organisations with
valuable data are likely to have such practices implemented across
a complex network of computers. In order to target such organisations,
a malicious actor must invest significant effort in circumventing
their security and spreading their ransomware across the computer
network. This typically involves multiple attack vectors, such as
targeted phishing, remote desktop protocol (RDP), and searching for
weak passwords\cite{stahie2020ransomware}. The malicious actor may
spend days or weeks escalating their level of access to the target's
network; in 2019, the mean dwell time (the duration a threat is present
in a system before it is detected) of ransomware was 43 days \cite{infocyte2019report}.
Only when the malicious actor has a level of access high enough to
compromise the target's backups will they start encrypting files,
ensuring that the only way to restore the data is by paying the ransom.
While this process of circumventing security requires a large investment
of time and effort, it enables the malicious actor to extort organisations
for ransoms far larger than previously possible \cite{coalition2020cyber}.

A key feature in any ransomware negotiation is the reliability of
the malicious actor; the victim's willingness to pay the ransom demand
must be affected by the likelihood of getting their data back afterwards.
There are two major reasons for a victim to not get their data after
paying. The first is that the malicious actor chooses not to return
the data, perhaps to avoid the cost incurred by doing so. This is
not a sound business practice, as it results in a loss of perceived
reliability that reduces a victim's willingness to pay, and hence,
the malicious actor's profits \cite{cartwright2019ransomware}, which
does not fit with the business-like stance that malicious actors operating
ransomware have adopted \cite{asokan2020experts}. Therefore, we assume
that our malicious actor will always attempt to return the victim's
data. But how does one unintentionally not return data? Modern ransomware
operates by encrypting the victim's data using public-private key
encryption, rendering the data unreadable until it has been decrypted
\cite{young1996cryptovirology}. In order to decrypt their data, the
victim must pay the malicious actor in exchange for the decryption
key. As the encryption phase of the attack must proceed rapidly to
avoid notice, flaws may arise during the process. If such flaws occur,
then the decryption key will fail to decrypt the affected data. This
results in some or all of the victim's files becoming permanently
irretrievable, which is unlikely to be discovered until after the
victim has paid the ransom. If a strain of ransomware is known to
have a history of failure, the target's willingness to pay is reduced,
and so the reliability of the malicious actor's ransomware is an important
factor in the negotiations. Ransomware strains can vary quite significantly
in their reliability, depending on how heavily the malicious actor
invested in the development of the ransomware. They may even invest
in providing ``customer service'' to their victims, walking them
through the decryption process to further improve their image of reliability
\cite{ng2017malware}.

The potential for negotiation adds a significant feature from game
theory to the dynamics of targeted ransomware; information asymmetry.
The two parties to the negotiation have different degrees of information
about each other, information that is crucial to determining the outcome
of the negotiation. The malicious actor does not know exactly how
much the target's data is worth to them. This is a significant factor,
as the malicious actor seeks to set the ransom as high as possible
in order to justify the significant investment they make in their
attack. Therefore, estimating the value of the target's files accurately
is of great importance. To aid them in this, the malicious actor can
make use of both publicly available data, including shareholder's
reports and valuations, and private data found on the target's computer
network, such as up-to-date finances and business plans. However,
uncovering and interpreting an organisation's private documents is
not effortless, and so producing an accurate estimate for the value
of the target's data requires further investment from the malicious
actor. While investing greater effort can lead to a more accurate
estimate, it is highly unlikely that the malicious actor will ever
achieve perfect accuracy. Even so, we might expect the target to be
at a distinct disadvantage in terms of information asymmetry. However,
as a result of developments in the computer security industry in response
to targeted ransomware, this is not entirely true. The potential for
negotiation of large ransoms has led to the emergence of organisations
that offer professional ransomware negotiation services \cite{rundle2020ransomware}.
The negotiators have experience in dealing with malicious actors behind
the various strains of ransomware, allowing them to negotiate effectively
on the behalf of targeted organisations. The negotiators also offer
specialised knowledge. This can include statistics such as the reliability
of a given strain of malware, but it can also include less quantifiable
information, such as how aggressively the malicious actor behind a
particular strain will negotiate. Operators of targeted ransomware
are typically willing to reduce their demand in the interest of getting
paid, but not all are equally receptive to negotiations. Some malicious
actors are likely to perceive a low counteroffer to their demand as
an insult, and may react aggressively to punish the target. The severity
of the reaction depends on how aggressive the malicious actor is,
which varies depending on individual and cultural factors. In extreme
cases, they may even react by abandoning the negotiations, and their
ransom, so that future targets will be less likely to negotiate. This
aggressive behaviour is a double-edged sword; if the malicious actor
is not aggressive enough, then their targets will not pay them a large
ransom. If they are too aggressive, their inclination to punish their
targets for perceived insults will cost them ransoms. In order to
be successful, the malicious actor must balance their aggression with
their investments in their strain of targeted ransomware. In the next
section, we construct a model of ransomware negotiation that is based
on these key features.

\section{Modelling}

We propose to study the dynamics of targeted ransomware by modelling
the negotiations as a two-player game. This approach was inspired
by Selten's analysis of a two-player game modelling the interaction
between a hostage taker and a hostage negotiator \cite{selten1988simple}.
Our two players are the attacker $A$ and the defender $D$. Player
$A$ is a malicious actor (or group of malicious actors) operating
a strain of targeted ransomware. Player $D$ is an organisation targeted
by player $A$, assisted by a professional negotiator hired to negotiate
the ransom.

\subsection{Player $A$'s investment}

As noted in the previous section, there are three areas in which player
$A$ must invest in order to pull off a successful targeted attack:
\begin{itemize}
\item The circumvention of player $D$'s security
\item The reliability of player $A$'s ransomware
\item The estimation of the value of player $D$'s data
\end{itemize}
We will not be considering player $A$'s investment in circumventing
player $D$'s security here. While it is likely to be an important
factor in the overall dynamics of the system, it has little bearing
on the negotiations. Once player $D$'s computer network has been
infected, this investment only affects player $A$'s net profit, and
affects player $D$ not at all. However, the other two investments
are very significant to the negotiations. 

\begin{figure}
\begin{centering}
\includegraphics[width=0.5\textwidth]{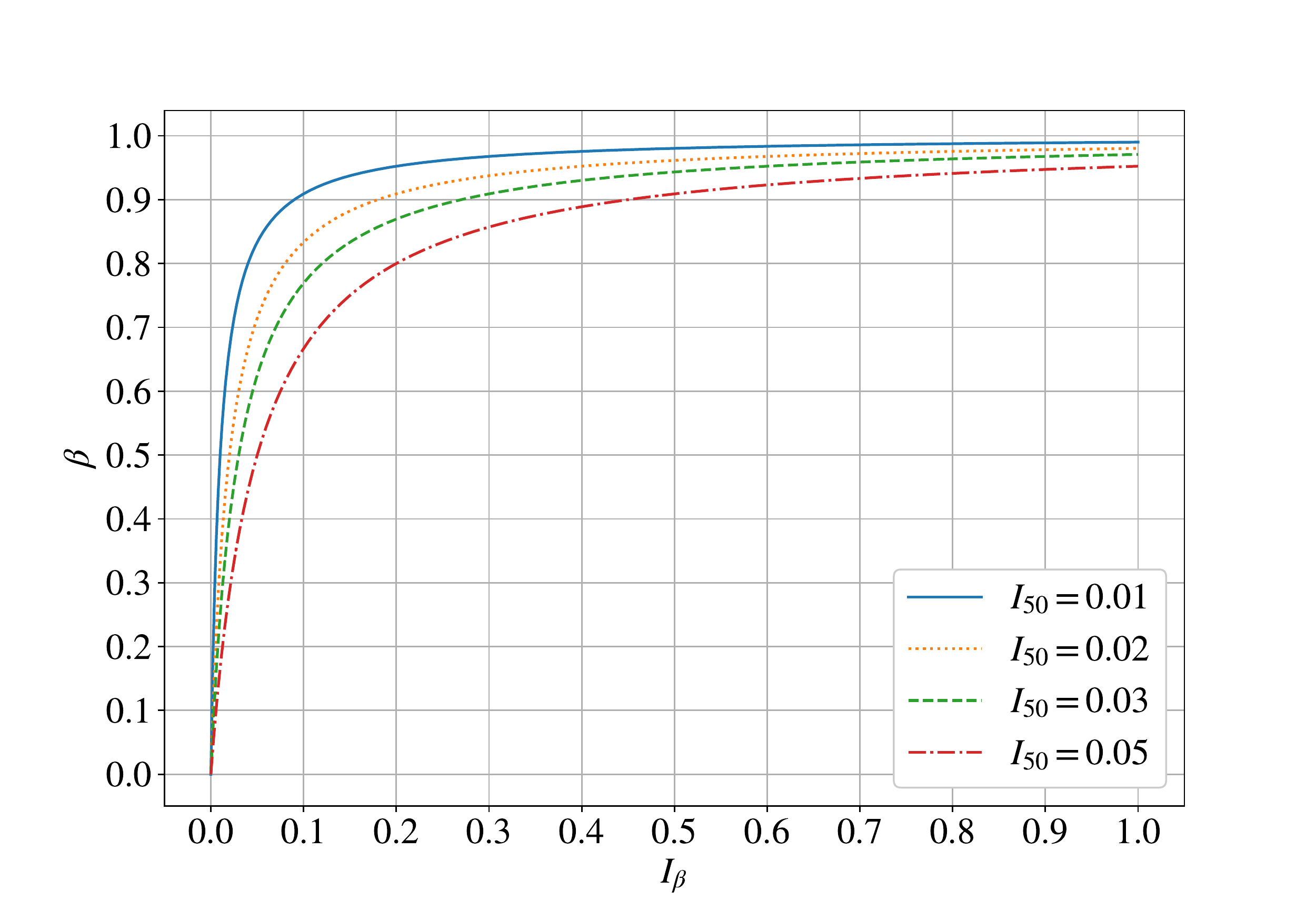}
\par\end{centering}
\caption{\label{fig:beta} $\beta$, the probability of succesful decryption,
depends on $I_{\beta}$ and $I_{50}$.}
\end{figure}

By the time player $A$ launches their attack, they have already developed
their strain of ransomware. In particular, they have invested in the
reliability of their ransomware. Investing in reliability is important;
if the decryption process is likely to fail, then player $D$ will
not be willing to pay very much for the decryption key. This is a
significant up-front development cost; it does not scale with the
number of targets player $A$ attacks. For the sake of simplicity,
we assume that player $A$ can amortize this investment over the targets
that they will infect with their strain of ransomware. We refer to
this investment as $I_{\beta}$. As $I_{\beta}$ increases, so does
the reliability of player $A$'s ransomware. Let $\beta$ be the probability
that the decryption key successfully decrypts data encrypted by the
ransomware. We choose $\beta$ such that
\begin{equation}
\beta=\frac{I_{\beta}}{I_{\beta}+I_{50}}\label{eq:beta}
\end{equation}
where $I_{50}$ is the amount of investment required to achieve a
reliability of $50\%$, or $\beta=0.5$. $I_{50}$ can be considered
as an economic scaling factor, and determines the amount of development
that player $A$ can purchase for a given investment. This choice
of $\beta$ results in diminishing return for large $I_{\beta}$,
so that $\beta\rightarrow1$ as $I_{\beta}\rightarrow\infty$, as
shown in Fig. \ref{fig:beta}. Note that $I_{50}$ and $I_{\beta}$
are dimensionless, so that our model remains generally applicable.
$I_{50}=0.02$ means that an investment of 2\% of the value of the
targeted files will yield a decryptor that is 50\% reliable. While
player $D$ doesn't know the value of $I_{\beta}$, their negotiator
can provide an estimate of $\beta$ from their experience of individual
ransomware strains.

\begin{figure}
\begin{centering}
\includegraphics[width=0.5\textwidth]{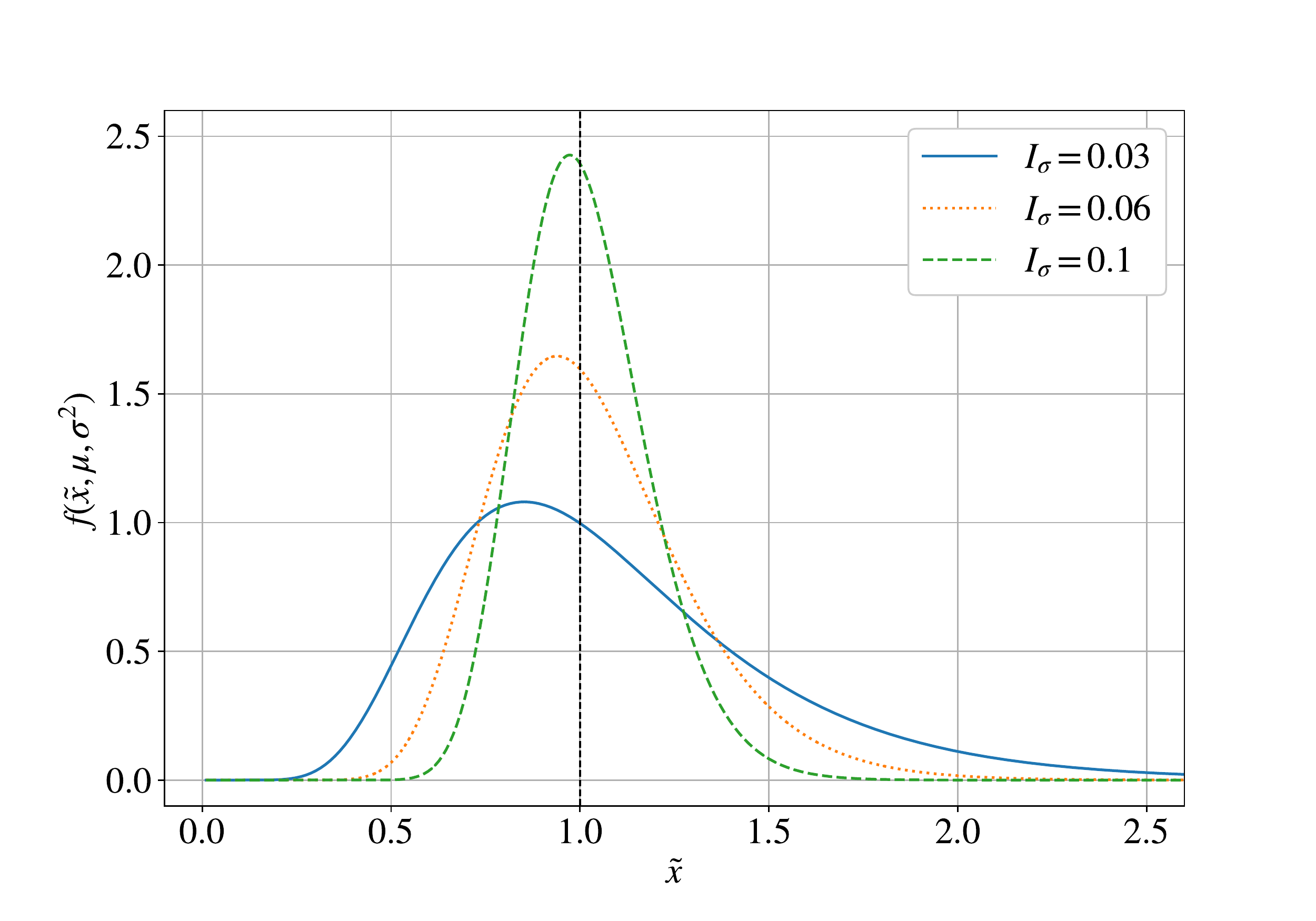}
\par\end{centering}
\caption{\label{fig:lognormal} Probability density function $f\left(\tilde{x},\mu,\sigma\right)$
for $x=1$, $I_{50}=0.02$ and varying levels of investment $I_{\sigma}$.}
\end{figure}

The final area in which player $A$ can invest is in their estimation
of how much player $D$'s data is worth. Without an accurate estimate,
player $A$ is unlikely to choose an optimal ransom demand, and so
investing in producing an accurate estimate is an important aspect
of their strategy. We let $x$ be the value that player $D$ attaches
to their encrypted data, and let $\tilde{x}$ be player $A$'s estimate
of $x$. We refer to player $A$'s investment in data value estimation
as $I_{\sigma}$. As $I_{\sigma}$ increases, so does the probability
that $\tilde{x}$ will be close to $x$. In this paper, we choose
to model $\tilde{x}$ as random variable following a $\text{Lognormal\ensuremath{\left(\mu,\sigma^{2}\right)}}$
distribution \cite{crow1987lognormal} with probability density function
\begin{equation}
f\left(\tilde{x},\mu,\sigma\right)=\frac{1}{\tilde{x}\sigma\sqrt{2\pi}}\exp\left(-\frac{\left(\ln\tilde{x}-\mu\right)^{2}}{2\sigma^{2}}\right)\label{eq:lognormal_pdf}
\end{equation}
 for $\tilde{x}>0$ and parameters

\begin{align}
\mu & =\ln x\nonumber \\
\sigma & =1-\frac{I_{\sigma}}{I_{50}+I_{\sigma}}\label{eq:mu_sigma}
\end{align}
With the chosen parameters, $\tilde{x}$ has median $x$. As $I_{\sigma}\rightarrow\infty$
, $\sigma\rightarrow0$, the variance $\left(e^{\sigma^{2}}-1\right)\left(e^{2\mu+\sigma^{2}}\right)\rightarrow0$
and mean $xe^{\frac{\sigma^{2}}{2}}\rightarrow x$. The Lognormal
distribution has previously been used for modelling positive quantities
that are determined by human behaviour \cite{gualandi2019human}.
As $\tilde{x}$ is player $A$'s estimate of how much player $D$
values their data, we believe that this is an appropriate choice.
Fig.  \ref{fig:lognormal} shows the probability density function
$f\left(\tilde{x},\mu,\sigma\right)$ for $x=1$ and varying levels
of investment $I_{\sigma}$. $I_{\sigma}$ is scaled similarly to
$I_{\beta}$. As $I_{\sigma}$ increases, the distribution narrows
around $x$.

\subsection{Negotiation}

Once player $A$ has infected player $D$'s computer network with
ransomware, the negotiation begins. Player $A$ issues a ransom demand
$R$, and player $D$ must decide how to respond. We now examine the
negotiation process between player $A$ and player $D$ that occurs
if player $D$ is willing to pay, but does not wish to pay the full
demand. This is a time-sensitive issue, particularly for player $D$.
Player $D$ cannot conduct business while their data is encrypted,
but they still incur costs, which can grow significant over a protracted
negotiation. Player $A$ does not suffer such ongoing costs, but they
have invested significant resources in the attack, and the longer
the negotiation continues, the more time player $D$ has to consider
their position. This makes a rapid negotiation process highly desirable,
and so we make a simplification to the negotiation process and model
it in the simplest way possible in the manner used by Selten \cite{selten1988simple}.
Player $A$ issues a ransom demand, player $D$ responds with a counteroffer
$C$, then player $A$ decides whether to accept $C$ and hand over
the decryption key, or reject $C$ and abandon the negotiations.

This is a substantially simplified description of the negotiation
process and should not be taken literally. In reality, there may be
a series of offers and counteroffers that take place over time. However,
player $D$ has finite capital with which they can absorb the costs
incurred by not being able to conduct business; a drawn out negotiation
for a lower ransom may be more costly than a prompt negotiation for
a higher ransom.

\begin{figure}
\begin{centering}
\includegraphics[width=0.5\textwidth]{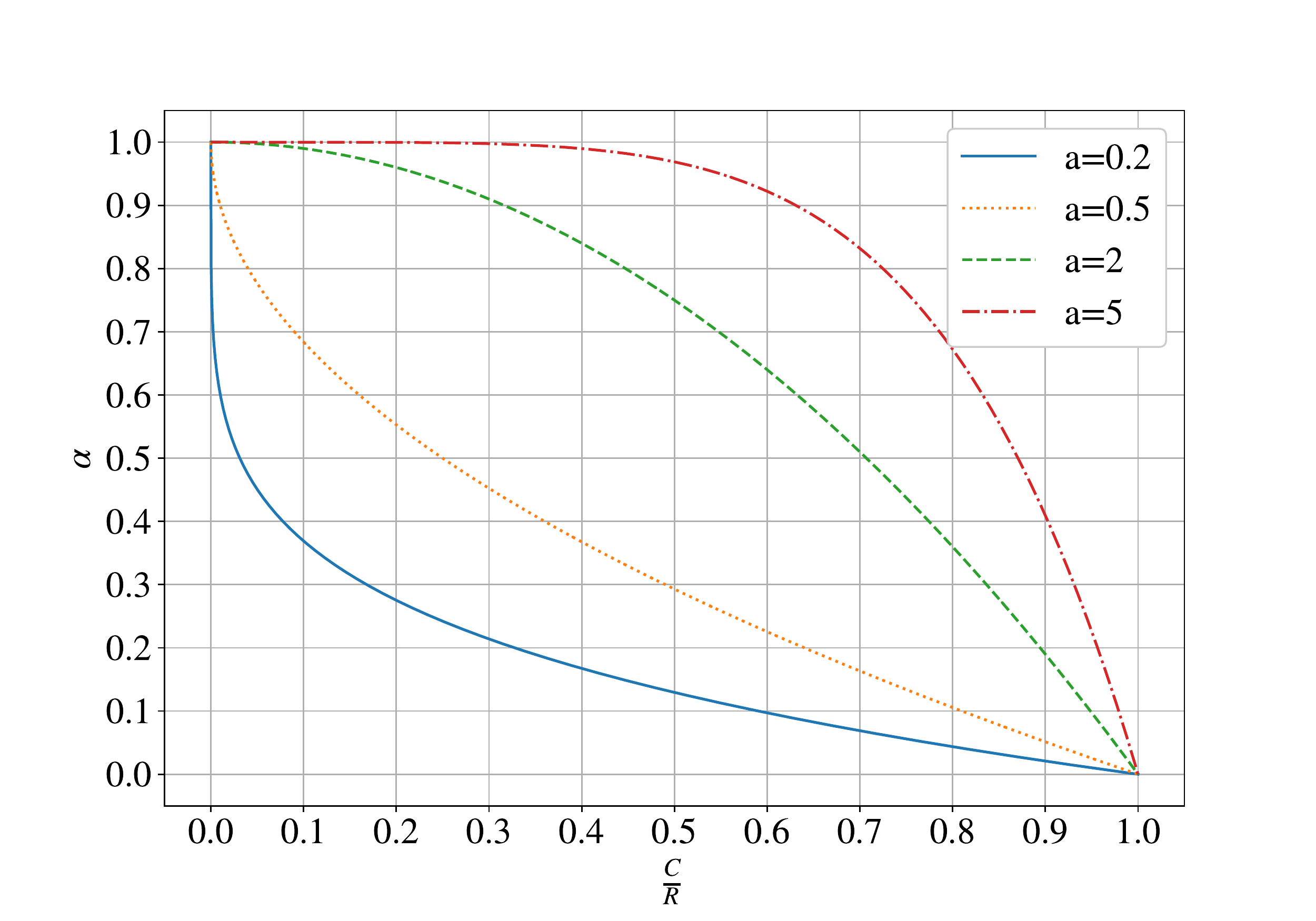}
\par\end{centering}
\caption{\label{fig:alpha}$\alpha$, the probability of an aggressive reaction
from player $A$, depends on $a$ and the ratio of counteroffer to
ransom demand $\frac{C}{R}$.}
\end{figure}

Why would player $A$ ever reject $C$? In doing so, they lose both
their potential earnings and waste any investment they've made in
the attack, which is clearly an undesirable outcome. However, player
$A$ must maintain their status as a threat; if they appear to be
willing to accept low counteroffers, they will only receive low counteroffers.
In order to maintain their credibility and their profits, player $A$
may punish player $D$ for making a low counteroffer. Therefore, we
must expect that with a positive probability $\alpha$, $A$ will
perceive a counteroffer $C<R$ as an affront, to which they react
aggressively by abandoning the negotiations. It is reasonable to suppose
that $\alpha$ will be greatest for $C=0$ and lowest for $C=R$.
We define $\alpha$ as
\begin{equation}
\alpha=1-\left(\frac{C}{R}\right)^{a}\label{eq:alpha}
\end{equation}
where $a>0$ is the aggression parameter of player $A$, quantifying
their tendency to perceive a low counteroffer as an affront and react
aggressively. This aggressive reaction is different to that implemented
by Selten \cite{selten1988simple}, where $\alpha=a\left(1-\frac{C}{D}\right)$
and $a\in\left[0,1\right]$ so that $\alpha\leq a\leq1$. The reason
for this is that in Selten's game, the aggressive reaction of the
hostage taker is to kill the hostage, while in our game, the aggressive
reaction is merely to not decrypt data. It is reasonable to assume
that a hostage taker, faced with having their ransom demand being
disregarded, might still refrain from killing their hostage. However,
a malicious actor, divorced from the consequences of their actions,
would have no reason not to react aggressively and abandon the negotiation.
This choice of $\alpha$ allows for a wide range of behaviour from
player $A$, with very lenient negotiations for $a<1$, scaling up
to very aggressive negotiations as $a$ increases. Larger $a$ causes
$\alpha$ to increase more rapidly as the difference between $C$
and $R$ increases, as shown in Fig. \ref{fig:alpha}. As with $\beta$,
player $D$ can estimate $a$ through the negotiator's experience
of interacting with individual ransomware operators.

If player $A$ does not react aggressively to a counteroffer $C<R$,
player $A$ receives payment $C$ and player $D$ receives the decryption
key, which they then attempt to use to decrypt their data. The probability
of successful decryption is $\beta$. If the decryption process is
successful, player $D$ retrieves their data. If the decryption is
unsuccessful, their data is rendered permanently irretrievable. The
possible outcomes of the game are detailed in Table \ref{tab:outcomes}.
In the next section, we consider the how these outcomes may arise
through rational decisions made by both players.

\begin{table}
\noindent \centering{}\caption{Payoff table for the targeted ransomware negotiation game.\label{tab:outcomes}}
\begin{tabular}{lcc}
\hline 
Outcome & \multicolumn{2}{c}{Payoff}\tabularnewline
\cline{2-3} \cline{3-3} 
 & Player $A$ & Player $D$\tabularnewline
\hline 
Aggressive rejection & $-I_{\beta}-I_{\sigma}$ & $-x$\tabularnewline
Decryption successful & $C-I_{\beta}-I_{\sigma}$ & $-C$\tabularnewline
Decryption failed & $C-I_{\beta}-I_{\sigma}$ & $-x-C$\tabularnewline
\hline 
\end{tabular}
\end{table}

\subsection{Summary of game rules}

The model variables are summarised in Table \ref{tab:variables}.
\begin{table*}
\noindent \centering{}\caption{Table of game variables with detail on the information asymmetry in
the targeted ransomware negotiation game. \label{tab:variables}}
\begin{tabular}{lccc}
\hline 
Variable & Description & Known to player $A$? & Known to player $D$?\tabularnewline
\hline 
$x$ & True value of $D$'s files & No & Yes\tabularnewline
$\tilde{x}$ & $A$'s estimate of $x$ & Yes & No\tabularnewline
$R$ & $A$'s ransom demand & Yes & Yes\tabularnewline
$C$ & $D$'s counteroffer & Yes & Yes\tabularnewline
$a$ & $A$'s aggression & Yes & Yes\tabularnewline
$I_{\beta}$ & $A$'s investment in reliability & Yes & No\tabularnewline
$I_{\sigma}$ & \multirow{1}{*}{$A$'s investment in estimating $x$} & Yes & No\tabularnewline
$\alpha$ & Probability that $A$ will react aggressively to $C<R$ & Yes & Yes\tabularnewline
$\beta$ & Probability that the decryption key works & Yes & Yes\tabularnewline
$\sigma$ & Scale parameter of distribution of $\tilde{x}$ & Yes & No\tabularnewline
\hline 
\end{tabular}
\end{table*}
The rules of the game are summarised as follows:
\begin{enumerate}
\item Player $A$ incurs cost $I_{\beta}+I_{\sigma}$ to infect player $D$'s
computer system and make a ransom demand $R$.
\item Player $D$ makes a counteroffer $C$.
\item Player $A$ aggressively rejects player $D$'s counteroffer with probability
$\alpha$.
\item If player $A$ does not aggressively reject the counteroffer, player
$A$ receives the counteroffer $C$ and player $D$ receives the decryption
key.
\item Player $D$'s data is successfully decrypted with probability $\beta$.
\end{enumerate}

\section{Analysis}

\subsection{Optimal choice of $C$}

\begin{figure}
\begin{centering}
\includegraphics[width=0.5\textwidth]{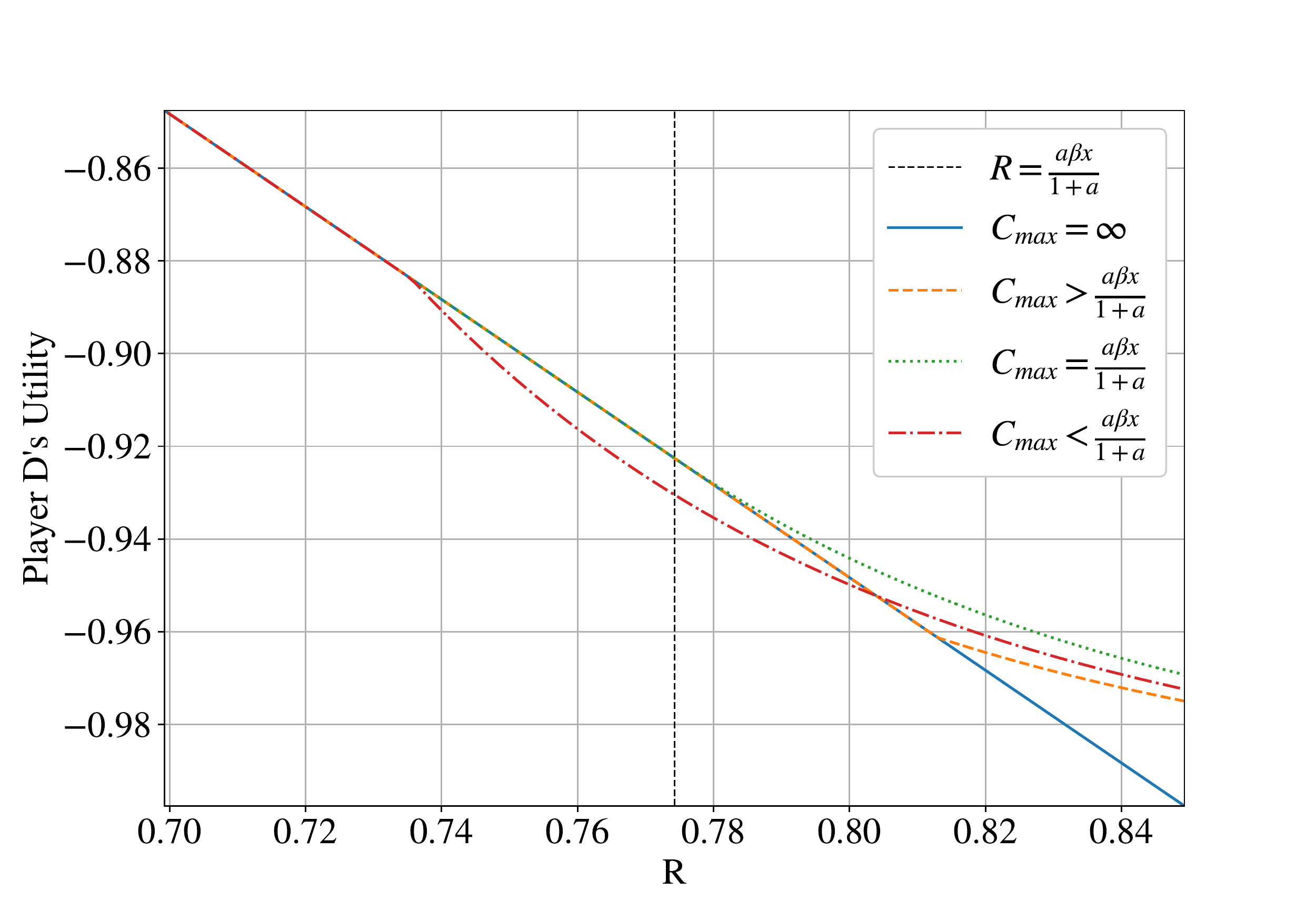}
\par\end{centering}
\caption{\label{fig:utility}Player $D$'s expected utility for varying $R$
and different values of $C_{max}$ where $a=10$, $I_{50}=0.02$ and
$I_{\beta}=0.1$.}
\end{figure}

In the subgame beginning with player $D$'s choice of $C$, player
$D$ knows that making a counteroffer $C<R$ will provoke an aggressive
reaction from player $A$ with probability $\alpha$. Player $D$
has no incentive to make a counteroffer $C>R$, so they rationally
chooses $C$ to maximize their expected utility $U$:
\[
U=\begin{cases}
-\left[R+\left(1-\beta\right)x\right] & C=R\\
-\left(1-\alpha\right)\left[C+\left(1-\beta\right)x\right]-\alpha x & C<R
\end{cases}
\]
The probability of successful decryption is $\beta$; therefore, the
expected value of the encrypted data is $\beta x$, and so player
$D$ will not offer more than $\beta x$ for the decryption key. We
assume that there exists a critical value $C_{max}\in\left[0,\beta x\right]$
such that player $D$'s expected utility is maximized by making counteroffer
$C_{max}$ if $R>C_{max}$. By substituting from Eq. (\ref{eq:alpha}),
player $D$'s expected utility for $C<R$ simplifies to
\[
U=-\left(\frac{C}{R}\right)^{a}\left[C-\beta x\right]-x
\]
To calculate the optimal value of $C$, we calculate
\[
\frac{\partial U}{\partial C}=\frac{C^{a-1}}{R^{a}}\left(a\beta x-C\left(1+a\right)\right)
\]
By setting $\frac{\partial U}{\partial C}=0$ we find that for $0<C<R$,
$U$ achieves its maximum at $C=\frac{a\beta x}{1+a}$. $\frac{\partial U}{\partial C}<0$
for $C>\frac{a\beta x}{1+a}$; for any $R>\frac{a\beta x}{1+a}$,
player $D$'s optimal counteroffer is $C=\frac{a\beta x}{1+a}$. $\frac{\partial U}{\partial C}>0$
for $C<\frac{a\beta x}{1+a}$; however, player $D$ will never make
a counteroffer $C>R$. If $R<\frac{a\beta x}{1+a}$, player $D$'s
optimal counteroffer is $C=R$. Thus, $C_{max}=\frac{a\beta x}{1+a}$,
yielding player $D$'s optimal counteroffer $\hat{C}$:
\begin{equation}
\hat{C}=\begin{cases}
R & R\leq\frac{a\beta x}{1+a}\\
\frac{a\beta x}{1+a} & R>\frac{a\beta x}{1+a}
\end{cases}\label{eq:c_hat}
\end{equation}
This scenario is demonstrated graphically in Fig. \ref{fig:utility}.
Any choice of $C_{max}\neq\frac{a\beta x}{1+a}$ results in a drop
in player $D$'s expected utility.

\subsection{Optimal choice of $R$}

In the subgame beginning with player $A$'s choice of $R$, player
$A$ knows that under rational decision-making, player $D$ will optimally
make counteroffer $\hat{C}$. Player $A$ rationally chooses $R$
to maximize their expected profit $P$:
\[
P=\begin{cases}
R-I_{\beta}-I_{\sigma} & R\leq C\\
\left(1-\alpha\right)C-I_{\beta}-I_{\sigma} & R>C
\end{cases}
\]
Substituting from Eq. (\ref{eq:c_hat}) yields
\begin{equation}
P=\begin{cases}
R-I_{\beta}-I_{\sigma} & R\leq\frac{a\beta x}{1+a}\\
\left(\frac{\frac{a\beta x}{1+a}}{R}\right)^{a}\left(\frac{a\beta x}{1+a}\right)-I_{\beta}-I_{\sigma} & R>\frac{a\beta x}{1+a}
\end{cases}\label{eq:p1}
\end{equation}
By differentiating with respect to $R$ we find that $\frac{\partial P}{\partial R}>0$
for $R<\frac{a\beta x}{1+a}$ and $\frac{\partial P}{\partial R}<0$
for $R>\frac{a\beta x}{1+a}$. Therefore, the optimal ransom demand
$\hat{R}=\frac{a\beta x}{1+a}$ is the highest ransom that player
$D$ is willing to pay. If player $A$ can reliably make demand $\hat{R}$,
player $D$ will always pay. Player $A$ will always make their maximum
profit, and there is no risk of an aggressive reaction from player
$A$, trivialising the negotiation. Under such conditions, player
$A$'s profit is $\frac{a\beta x}{1+a}-I_{\beta}-I_{\sigma}$. This
would suggest that, in order to maximise their profit, player $A$
should be infinitely aggressive, rejecting any counteroffer even slightly
lower than their demand (i.e. $a\rightarrow\infty$). 

\begin{figure}[h]
\begin{centering}
\includegraphics[width=0.5\textwidth]{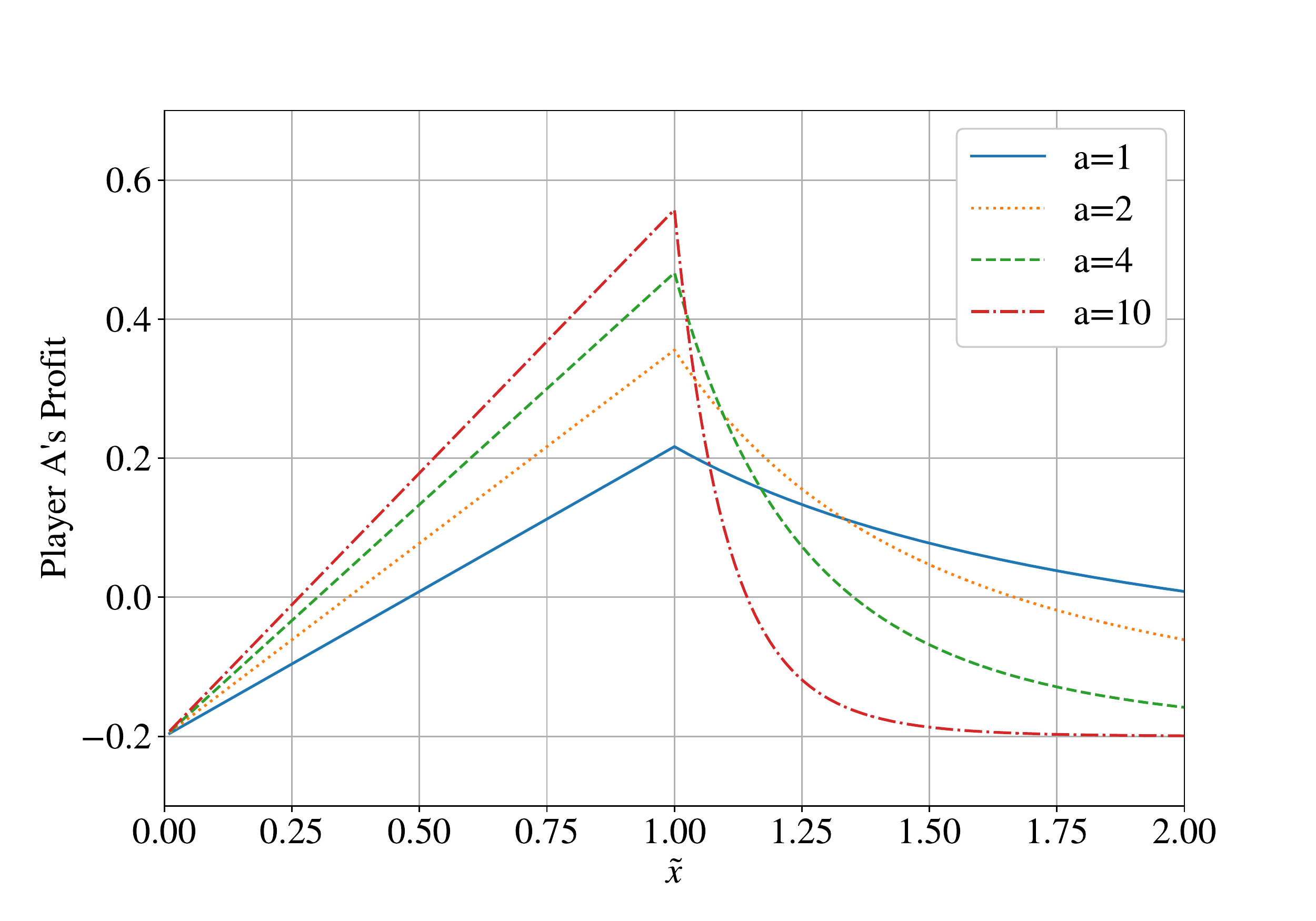}
\par\end{centering}
\caption{Player $A$'s expected profit as a function of $\tilde{x}$ for varying
$a$ when $x=1$, $I_{50}=0.02$, and $I_{\beta}=I_{\sigma}=0.1$.\label{fig:Attacker's-net-profit}}
\end{figure}

Of course, this ``ideal'' scenario is unrealistic, as it ignores
the often-significant effect of imperfect information \cite{kreps1982reputation,barrachina2014entry,durkota2015approximate}.
In reality, negotiations are not trivial affairs, and the risk of
an aggressive reaction is always present. Player $A$ does not know
$x$, only $\tilde{x}$, which, lacking any alternative, is what they
use to calculate their ransom demand. Therefore, under optimal play
while accounting for imperfect information, player $A$'s ransom demand
is $R=\frac{a\beta\tilde{x}}{1+a}$. The potential for error in player
$A$'s estimate $\tilde{x}$ gives rise to the necessity of negotiations
that may result in an aggressive reaction. We can illustrate this
by substituting $R=\frac{a\beta\tilde{x}}{1+a}$ into Eq. (\ref{eq:p1}).
Under optimal play,
\[
P=\begin{cases}
\frac{a\beta\tilde{x}}{1+a}-I_{\beta}-I_{\sigma} & \frac{a\beta\tilde{x}}{1+a}\leq\frac{a\beta x}{1+a}\\
\left(\frac{\frac{a\beta x}{1+a}}{\frac{a\beta\tilde{x}}{1+a}}\right)^{a}\left(\frac{a\beta x}{1+a}\right)-I_{\beta}-I_{\sigma} & \frac{a\beta\tilde{x}}{1+a}>\frac{a\beta x}{1+a}
\end{cases}
\]
which factors to 
\begin{equation}
P=\frac{a\beta}{1+a}\left.\begin{cases}
\tilde{x} & \tilde{x}\leq x\\
x\left(\frac{x}{\tilde{x}}\right)^{a} & \tilde{x}>x
\end{cases}\right]-I_{\beta}-I_{\sigma}\label{eq:p2}
\end{equation}
Thus, it is the potential for error in the estimate $\tilde{x}$ which
prevents player $A$ from playing optimally. The effect of error in
$\tilde{x}$ on player $A$'s profit is shown in Fig.  \ref{fig:Attacker's-net-profit}.
Here, for investment levels are fixed at $I_{\beta}=I_{\sigma}=0.1$,
so that the maximum profit depends on $a$. The effect of error differs
depending on whether player $A$ underestimates or overestimates the
value of the data. If they underestimate $x$, then their profit decreases
linearly with $\tilde{x}$. If they overestimate $x$, then the possibility
of an aggressive reaction emerges, which increases with both $a$,
and the error in $\tilde{x}$. High aggression might increase player
$A$'s capacity for demanding large ransoms, but at an increased risk
of an aggressive reaction.
In order to further understand how error and aggression interact,
we consider player $A$'s expected profit. First, let player $A$'s
gross profit be

\[
GP\left(\tilde{x},x;a,\beta\right)=\frac{a\beta}{a+1}\begin{cases}
\tilde{x} & \tilde{x}\leq x\\
x\left(\frac{x}{\tilde{x}}\right)^{a} & \tilde{x}>x
\end{cases}
\]
We let $x$ follow a distribution with probability density function
$g\left(x,M\right)$ for $x>0$ and mean $M$; the shape of the distribution
is unimportant in this case. With probability density functions $g\left(x,M\right)$
for $x$ and $f\left(\tilde{x},\mu,\sigma\right)$ for $\tilde{x}$,
Player $A's$ net profit for a given combination of aggression and
investments can be written as a convolution of $x$ and $\tilde{x}$
in double integral form. However, due to our choice of $\alpha$ and
Lognormal $\tilde{x}$, by making the substitution $y=\frac{\tilde{x}}{x}$
we can reduce the double integral to a single integral. 

\begin{figure}
\begin{centering}
\includegraphics[width=0.5\textwidth]{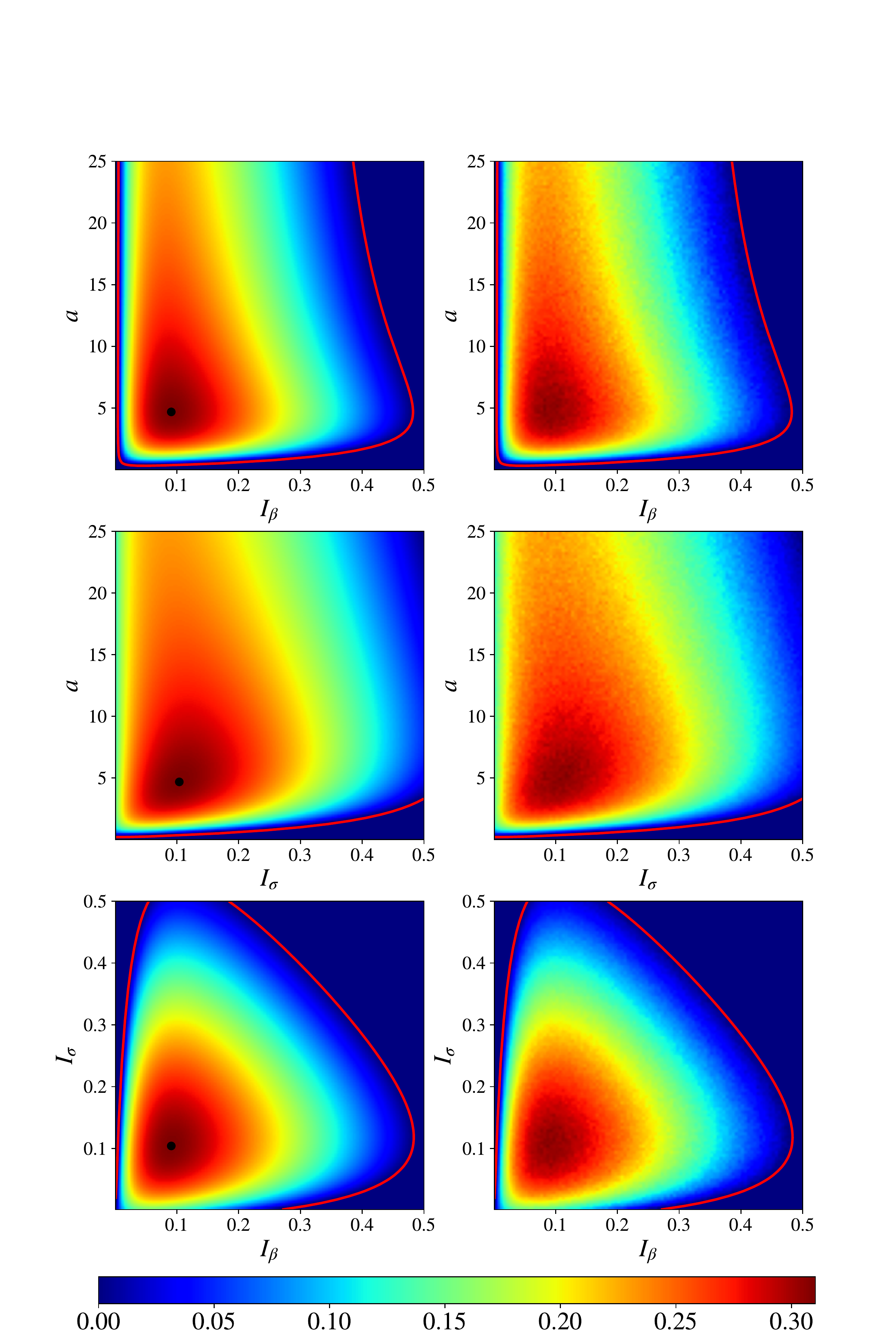}
\par\end{centering}
\caption{Player $A$'s profit for varying strategy parameters $\left(a,I_{\beta},I_{\sigma}\right)$
when $I_{50}=0.02$. In each plot, the hidden parameter is set to
its optimal value. The maximum mean profit achieved is marked by a
black dot. The red curve marks where the the mean profit is equal
to 0.\label{fig:profit_master}}
\end{figure}

\begin{widetext}\begin{align} P\left(a,I_{\beta},I_{\sigma}\right)= & \int_{x=0}^{x=\infty}\int_{y=0}^{y=\infty}GP\left(xy,x;a,\beta\right)f\left(xy,\mu,\sigma^{2}\right)g\left(x,M\right)xdydx-I_{\beta}-I_{\sigma}\nonumber \\ = & \int_{x=0}^{x=\infty}\int_{y=0}^{y=\infty}\left(\frac{a\beta}{a+1}\right)\left.\begin{cases} xy & xy\leq x\\ x\left(\frac{1}{y}\right)^{a} & xy>x \end{cases}\right]\frac{1}{xy\sigma\sqrt{2\pi}}e^{-\frac{\left[\ln\left(xy\right)-\ln\left(x\right)\right]^{2}}{2\sigma^{2}}}g\left(x,M\right)xdydx-I_{\beta}-I_{\sigma}\nonumber \\ = & \left(\frac{a\beta}{a+1}\right)\int_{x=0}^{x=\infty}g\left(x,M\right)xdx\int_{y=0}^{y=\infty}\left.\begin{cases} y & y\leq1\\ y^{-a} & y>1 \end{cases}\right]\frac{1}{y\sigma\sqrt{2\pi}}e^{-\frac{\ln^{2}y}{2\sigma^{2}}}dy-I_{\beta}-I_{\sigma}\nonumber \\ = & \left(\frac{a\beta}{a+1}\right)M\left[\int_{y=0}^{y=1}y\frac{1}{y\sigma\sqrt{2\pi}}e^{-\frac{\ln^{2}y}{2\sigma^{2}}}dy+\int_{y=1}^{y=\infty}y^{-a}\frac{1}{y\sigma\sqrt{2\pi}}e^{-\frac{\ln^{2}y}{2\sigma^{2}}}dy\right]-I_{\beta}-I_{\sigma}\label{eq:profit_master} \end{align}\end{widetext}By
making this substitution, we can see that what appears to be a convolution
of $x$ and $\tilde{x}$ depends merely on the ratio $\frac{\tilde{x}}{x}$.
Player $A$'s strategy is consistent across all values of $x$, so
it is only the mean $M$ of the distribution of $x$ that remains
in the final expression. In this form, we can clearly see where each
element of player $A$'s strategy $\left(a,I_{\beta},I_{\sigma}\right)$
comes into play. As aggression $a$ increases, the multiplicative
term $\frac{a}{a+1}$ increases, but the second integral in the sum,
where player $A$ has overestimated $x$, converges to $0$. Increasing
$I_{\beta}$ increases costs, but leads to increased $\beta$ which
may increase profit. Increasing $I_{\sigma}$ also increases costs,
but narrows the distribution of $\tilde{x}$ around $x$, allowing
for greater aggression at decreased risk of aggressive reaction. Thus,
through our choice of $\alpha$ and $\tilde{x}$, we can more clearly
demonstrate how player $A$'s strategy depends on the interaction
between the various elements of their strategy. The optimal counteroffer
and expected profit for parameters corresponding to various strategies
are shown in Table  \ref{tab:strategies}. Fig.  \ref{fig:profit_master}
shows player $A$'s profit for varying parameters $a$, $I_{\beta}$
and $I_{\sigma}$. The left column shows results calculated via numerical
integration of Eq. (\ref{eq:profit_master}) for $M=1$. The right
column verifies these calculated results with an agent-based simulation,
where results are averaged over $n=10000$ runs with constant target
data value $x=1$. These figures demonstrate that player $A$'s optimal
strategy is $\left(a,I_{\beta},I_{\sigma}\right)=\left(4.68,0.091,0.104\right)$,
rather than the naive maximal aggression $a\rightarrow\infty$ noted
previously; in order to realise their potential profits, the attacker
must be willing to negotiate.

\begin{table}
\begin{centering}
\caption{Model variables for optimal and various sub-optimal strategies.\label{tab:strategies}}
\begin{tabular}{|c|c|c|c|}
\hline 
Strategy Type & $\left(a,I_{\beta},I_{\sigma}\right)$ & $C$ & $P$\tabularnewline
\hline 
\hline 
Optimal & $\left(4.68,0.091,0.104\right)$ & $0.675$ & $0.304$\tabularnewline
\hline 
Low Aggression & $\left(2.34,0.091,0.104\right)$ & $0.574$ & $0.276$\tabularnewline
\hline 
High Aggression & $\left(9.36,0.091,0.104\right)$ & $0.741$ & $0.284$\tabularnewline
\hline 
Low Reliability & $\left(4.68,0.041,0.104\right)$ & $0.554$ & $0.265$\tabularnewline
\hline 
High Reliability & $\left(4.68,0.182,0.104\right)$ & $0.742$ & $0.264$\tabularnewline
\hline 
Low Accuracy & $\left(4.68,0.091,0.052\right)$ & $0.675$ & $0.283$\tabularnewline
\hline 
High Accuracy & $\left(4.68,0.091,0.208\right)$ & $0.675$ & $0.267$\tabularnewline
\hline 
\end{tabular}
\par\end{centering}

\end{table}

\section{Conclusion}

We have constructed a simple game-theory model of targeted ransomware
negotiations between two agents; an attacker, a malicious actor operating
targeted ransomware, and a defender, their target. Our model focuses
on three key elements of the attacker's strategy which have a direct
effect on the negotiations; aggression, investment in the reliability
of their ransomware, and investment in their estimation of target
data value. By thoroughly analysing our model, we demonstrate how
the necessity of negotiation arises from optimal decision made by
the agents under imperfect information. We show how the key elements
of the attacker's strategy interact with each other, and demonstrate
by numerical integration and by agent-based simulation that the attacker's
profit depends on developing a balance of investment and aggression.
While ransomware strains do show a significant variety of behaviour,
even within the subgroup of targeted ransomware, the features present
in this model are quite generic, and so we expect the insights provided
to be broadly applicable in the study of targeted ransomware. 

In this paper we have focused solely on the negotiations of targeted
ransomware that follow the infection of an organisation in order to
construct a simple model that can be analysed in detail. In this scenario,
the attacker has greater agency to determine the outcome, and so our
analysis is focused on the decisions of the attacker. As noted in
the Modelling section, this model neglects the attacker's investment
in circumventing security, as it has no direct effect on the negotiations.
By extension, the defender's investment in security is omitted. One
could grant greater agency to the defender by considering how the
defender's investment in security acts as a deterrent to a potential
attacker. Every additional effort, such as redundant backups, data
exfiltration countermeasures, and employee training makes an organisation
harder to hold to ransom, and so less attractive as a target. Such
deterrents, and other solutions to the growing problem of ransomware,
are of vital importance to strengthening our shared internet security.
Our model contributes to this goal by improving our understanding
of various negotiation strategies through the use of game theory.
As cybercrime continues to threaten our increasingly technology-dependent
lives, we hope that the application of game theory to targeted ransomware
will be of interest to a wide audience.

\section{Conflict of Interest Statement}

This research was supported by the Irish Research Council and McAfee
LLC through the Irish Research Council Employment-Based Postgraduate
Programme. Pierce Ryan is an employee of McAfee LLC. Any opinions
expressed in this work are solely those of the authors, and do not
necessarily reflect the views of the supporting organisations.

\bibliographystyle{IEEEtran}
\bibliography{ransomware_paper}

\begin{thebibliography}{10}
\providecommand{\url}[1]{#1}
\csname url@samestyle\endcsname
\providecommand{\newblock}{\relax}
\providecommand{\bibinfo}[2]{#2}
\providecommand{\BIBentrySTDinterwordspacing}{\spaceskip=0pt\relax}
\providecommand{\BIBentryALTinterwordstretchfactor}{4}
\providecommand{\BIBentryALTinterwordspacing}{\spaceskip=\fontdimen2\font plus
\BIBentryALTinterwordstretchfactor\fontdimen3\font minus
  \fontdimen4\font\relax}
\providecommand{\BIBforeignlanguage}[2]{{%
\expandafter\ifx\csname l@#1\endcsname\relax
\typeout{** WARNING: IEEEtran.bst: No hyphenation pattern has been}%
\typeout{** loaded for the language `#1'. Using the pattern for}%
\typeout{** the default language instead.}%
\else
\language=\csname l@#1\endcsname
\fi
#2}}
\providecommand{\BIBdecl}{\relax}
\BIBdecl

\bibitem{kalaimannan2017influences}
E.~Kalaimannan, S.~K. John, T.~DuBose, and A.~Pinto, ``Influences on
  ransomware's evolution and predictions for the future challenges,''
  \emph{Journal of Cyber Security Technology}, vol.~1, pp. 23--31, 2017.

\bibitem{maigida2019systematic}
A.~M. Maigida, S.~M. Abdulhamid, M.~Olalere \emph{et~al.}, ``Systematic
  literature review and metadata analysis of ransomware attacks and detection
  mechanisms,'' \emph{Journal of Reliable Intelligent Environments}, vol.~5,
  pp. 67--89, 2019.

\bibitem{beek2016targeted}
C.~Beek and A.~Furtak. (2016, February) Targeted ransomware no longer a future
  threat.
  www.mcafee.com/enterprise/en-us/assets/reports/rp-targeted-ransomware.pdf.
  Accessed 7 Oct 2020.

\bibitem{bajpai2020dissecting}
P.~Bajpai and R.~Enbody, ``Dissecting. net ransomware: key generation,
  encryption and operation,'' \emph{Network Security}, vol. 2020, no.~2, pp.
  8--14, 2020.

\bibitem{coveware2020ransomware}
Coveware. (2020, November) Ransomware demands continue to rise as data
  exfiltration becomes common, and maze subdues.
  www.coveware.com/blog/q3-2020-ransomware-marketplace-report. Accessed 25 Nov
  2020.

\bibitem{zimba2019economic}
A.~Zimba and M.~Chishimba, ``On the economic impact of crypto-ransomware
  attacks: the state of the art on enterprise systems,'' \emph{European Journal
  for Security Research}, vol.~4, no.~1, pp. 3--31, 2019.

\bibitem{myerson2013game}
R.~B. Myerson, \emph{Game theory}.\hskip 1em plus 0.5em minus 0.4em\relax
  Harvard university press, 2013.

\bibitem{brown1999ecology}
J.~S. Brown, J.~W. Laundr{\'e}, and M.~Gurung, ``The ecology of fear: optimal
  foraging, game theory, and trophic interactions,'' \emph{Journal of
  mammalogy}, vol.~80, pp. 385--399, 1999.

\bibitem{mcgill2007evolutionary}
B.~J. McGill and J.~S. Brown, ``Evolutionary game theory and adaptive dynamics
  of continuous traits,'' \emph{Annual Review of Ecology, Evolution, and
  Systematics}, vol.~38, pp. 403--435, 2007.

\bibitem{smith1973logic}
J.~M. Smith and G.~R. Price, ``The logic of animal conflict,'' \emph{Nature},
  vol. 246, pp. 15--18, 1973.

\bibitem{friedman1998economic}
D.~Friedman, ``On economic applications of evolutionary game theory,''
  \emph{Journal of evolutionary economics}, vol.~8, pp. 15--43, 1998.

\bibitem{selten1990bounded}
R.~Selten, ``Bounded rationality,'' \emph{Journal of Institutional and
  Theoretical Economics (JITE)/Zeitschrift f{\"u}r die gesamte
  Staatswissenschaft}, vol. 146, pp. 649--658, 1990.

\bibitem{lukas2012earnouts}
E.~Lukas, J.~J. Reuer, and A.~Welling, ``Earnouts in mergers and acquisitions:
  A game-theoretic option pricing approach,'' \emph{European Journal of
  Operational Research}, vol. 223, pp. 256--263, 2012.

\bibitem{cerdeiro2017contagion}
D.~A. Cerdeiro, ``Contagion exposure and protection technology,'' \emph{Games
  and Economic Behavior}, vol. 105, pp. 230--254, 2017.

\bibitem{kydd1997game}
A.~Kydd, ``Game theory and the spiral model,'' \emph{World Politics}, vol.~49,
  pp. 371--400, 1997.

\bibitem{ward1993game}
H.~Ward, ``Game theory and the politics of the global commons,'' \emph{Journal
  of Conflict Resolution}, vol.~37, pp. 203--235, 1993.

\bibitem{caulfield2015optimizing}
T.~Caulfield and A.~Fielder, ``Optimizing time allocation for network
  defence,'' \emph{Journal of Cybersecurity}, vol.~1, pp. 37--51, 2015.

\bibitem{lindsay2015tipping}
J.~R. Lindsay, ``Tipping the scales: the attribution problem and the
  feasibility of deterrence against cyberattack,'' \emph{Journal of
  Cybersecurity}, vol.~1, pp. 53--67, 2015.

\bibitem{laszka2017economics}
A.~Laszka, S.~Farhang, and J.~Grossklags, ``On the economics of ransomware,''
  in \emph{International Conference on Decision and Game Theory for
  Security}.\hskip 1em plus 0.5em minus 0.4em\relax Springer, 2017, pp.
  397--417.

\bibitem{cartwright2019pay}
E.~Cartwright, J.~Hernandez~Castro, and A.~Cartwright, ``To pay or not: game
  theoretic models of ransomware,'' \emph{Journal of Cybersecurity}, vol.~5, p.
  tyz009, 2019.

\bibitem{hu2020optimal}
H.~Hu, Y.~Liu, C.~Chen, H.~Zhang, and Y.~Liu, ``Optimal decision making
  approach for cyber security defense using evolutionary game,'' \emph{IEEE
  Transactions on Network and Service Management}, vol.~17, no.~3, pp.
  1683--1700, 2020.

\bibitem{caporusso2018game}
N.~Caporusso, S.~Chea, and R.~Abukhaled, ``A game-theoretical model of
  ransomware,'' in \emph{International Conference on Applied Human Factors and
  Ergonomics}.\hskip 1em plus 0.5em minus 0.4em\relax Springer, 2018, pp.
  69--78.

\bibitem{hernandez2020economic}
J.~Hernandez-Castro, A.~Cartwright, and E.~Cartwright, ``An economic analysis
  of ransomware and its welfare consequences,'' \emph{Royal Society Open
  Science}, vol.~7, p. 190023, 2020.

\bibitem{cartwright2019ransomware}
A.~Cartwright and E.~Cartwright, ``Ransomware and reputation,'' \emph{Games},
  vol.~10, p.~26, 2019.

\bibitem{li2021game}
Z.~Li and Q.~Liao, ``Game theory of data-selling ransomware,'' \emph{Journal of
  Cyber Security and Mobility}, pp. 65--96, 2021.

\bibitem{panda2019targeted}
P.~Security. (2019, July) Targeted ransomware attacks: the easy choice for
  cybercriminals.
  www.pandasecurity.com/en/mediacenter/security/targeted-ransomware/. Accessed
  4 Feb 2021.

\bibitem{stahie2020ransomware}
S.~Stahie. (2020, March) Why is targeted ransomware so dangerous?
  businessinsights.bitdefender.com/why-is-targeted-ransomware-so-dangerous.
  Accessed 14 Jan 2021.

\bibitem{infocyte2019report}
C.~Gerritz. (2019) 2019 mid-market threat and incident response report.
  www.infocyte.com/resources/mid-market-threat-and-incident-response-report/.
  Accessed 14 Jan 2020.

\bibitem{coalition2020cyber}
Coalition. (2020) Cyber insurance claims report.
  info.coalitioninc.com/download-2020-cyber-claims-report.html. Accessed 14 Jan
  2021.

\bibitem{asokan2020experts}
A.~Asokan. (2020, February) Experts warn: Targeted ransomware attacks to surge.
  www.bankinfosecurity.com/experts-warn-targeted-ransomware-attacks-to-surge-a-13782.
  Accessed 14 Jan 2021.

\bibitem{young1996cryptovirology}
A.~Young and M.~Yung, ``Cryptovirology: Extortion-based security threats and
  countermeasures,'' in \emph{Proceedings 1996 IEEE Symposium on Security and
  Privacy}.\hskip 1em plus 0.5em minus 0.4em\relax IEEE, 1996, pp. 129--140.

\bibitem{ng2017malware}
A.~Ng. (2017, July) Malware now comes with customer service.
  www.cnet.com/news/ransomware-goes-pro-customer-service-google-25-million-black-hat/.
  Accessed 14 Jan 2021.

\bibitem{rundle2020ransomware}
J.~Rundle. (2020, August) Ransomware negotiators gain prominence as attacks
  increase.
  www.wsj.com/articles/ransomware-negotiators-gain-prominence-as-attacks-increase-11598866201.
  Accessed 14 Jan 2021.

\bibitem{selten1988simple}
R.~Selten, ``A simple game model of kidnapping,'' in \emph{Models of strategic
  rationality}.\hskip 1em plus 0.5em minus 0.4em\relax Springer, 1988, pp.
  77--93.

\bibitem{crow1987lognormal}
E.~L. Crow and K.~Shimizu, \emph{Lognormal distributions}.\hskip 1em plus 0.5em
  minus 0.4em\relax Marcel Dekker New York, 1987.

\bibitem{gualandi2019human}
S.~Gualandi and G.~Toscani, ``Human behavior and lognormal distribution. a
  kinetic description,'' \emph{Mathematical Models and Methods in Applied
  Sciences}, vol.~29, pp. 717--753, 2019.

\bibitem{kreps1982reputation}
D.~M. Kreps and R.~Wilson, ``Reputation and imperfect information,''
  \emph{Journal of economic theory}, vol.~27, no.~2, pp. 253--279, 1982.

\bibitem{barrachina2014entry}
A.~Barrachina, Y.~Tauman, and A.~Urbano, ``Entry and espionage with noisy
  signals,'' \emph{Games and economic behavior}, vol.~83, pp. 127--146, 2014.

\bibitem{durkota2015approximate}
K.~Durkota, V.~Lis{\`y}, B.~Bo{\v{s}}ansk{\`y}, and C.~Kiekintveld,
  ``Approximate solutions for attack graph games with imperfect information,''
  in \emph{International Conference on Decision and Game Theory for
  Security}.\hskip 1em plus 0.5em minus 0.4em\relax Springer, 2015, pp.
  228--249.

\end{thebibliography}

\begin{IEEEbiography}[{\includegraphics[width=1in,height=1.25in,clip,keepaspectratio]{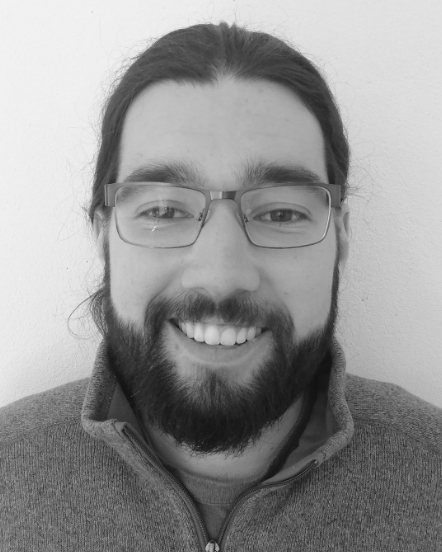}}]{Pierce Ryan} (Graduate Student Member, IEEE) received a BSc in mathematical sciences from University College Cork, Ireland. He is currently pursuing a PhD in applied mathematics at University College Cork in partnership with McAfee LLC under the Irish Research Council Employment-Based Postgraduate Programme while working as a Malware Researcher on the Artificial Intelligence Research Team at McAfee LLC. His research interests include data modelling and stochastic modelling for computer security systems, as well as dynamical systems incorporating time-delay, forcing and switches. \end{IEEEbiography}

\begin{IEEEbiography}[{\includegraphics[width=1in,height=1.25in,clip,keepaspectratio]{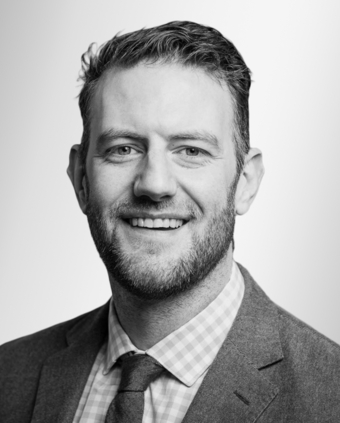}}]{John Fokker} received a BA in information and communication technology from Amsterdam University of Applied Sciences, Netherlands. He received a Masters in criminal investigation from Politieacademie Apeldoorn, Netherlands. He is currently a Principal Engineer \& Head of Cyber Investigations at Trellix Labs. His research interests include threat intelligence, malware analysis,  and the cyber criminal Underground. John is also one of the co-founders of the NoMoreRansom project. \end{IEEEbiography}

\begin{IEEEbiography}[{\includegraphics[width=1in,height=1.25in,clip,keepaspectratio]{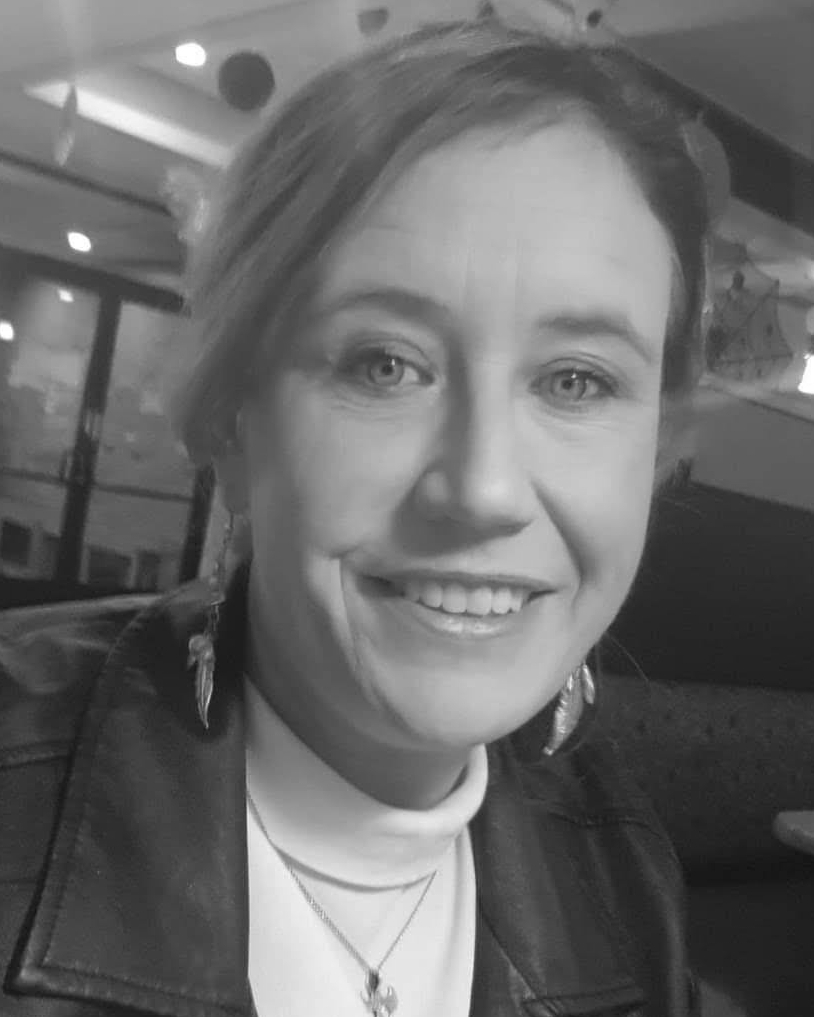}}]{Sorcha Healy} received a BSc in applied mathematics and physics from Dublin Institute of Technology, Ireland. She received a PhD in computational physics from University College Cork, Ireland. She is currently a Principal Data Scientist on the Global Data Science \& Analytics team at Microsoft Ireland. Her current research interests include natural language processing and customer product feedback. \end{IEEEbiography}

\begin{IEEEbiography}[{\includegraphics[width=1in,height=1.25in,clip,keepaspectratio]{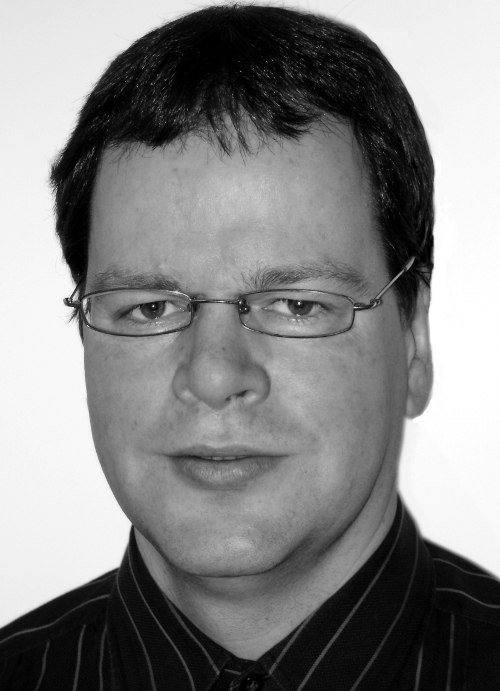}}]{Andreas Amann} received the Diploma degree in physics from the University of Bonn, Germany, and the PhD degree in physics from the Technical University of Berlin, Germany. After being a visiting researcher at the Institute of Applied Optics, Florence, Italy, and a postdoctoral researcher at Tyndall National Institute, Cork, Ireland, he is now a Senior Lecturer in the School of Mathematical Sciences at University College Cork, Ireland. His research interest include nonlinear dynamics, mathematical modelling and machine learning. \end{IEEEbiography}
\EOD 
\end{document}